\documentclass[titlepage,12pt]{article}

\usepackage{graphicx}
\usepackage{amssymb}

\newcommand{\R}{\mathbb R}
\newcommand{\C}{\mathbb C}

\begin{document}

\baselineskip=18pt

\begin{center}
{\Large{\bf Stability and Hopf Bifurcation in an \\ Hexagonal
Governor System With a Spring}}
\end{center}

\vspace{0.5cm}

\begin{center}
{\large Jorge Sotomayor}
\end{center}
\begin{center}
{\em Instituto de Matem\'atica e Estat\'{\i}stica, Universidade de
S\~ao Paulo\\ Rua do Mat\~ao 1010, Cidade Universit\'aria\\ CEP
05.508-090, S\~ao Paulo, SP, Brazil
\\}e--mail:sotp@ime.usp.br
\end{center}

\begin{center}
{\large Luis Fernando Mello}
\end{center}
\begin{center}
{\em Instituto de Ci\^encias Exatas, Universidade Federal de
Itajub\'a\\Avenida BPS 1303, Pinheirinho, CEP 37.500-903,
Itajub\'a, MG, Brazil
\\}e--mail:lfmelo@unifei.edu.br
\end{center}

\begin{center}
{\large Denis de Carvalho Braga}
\end{center}
\begin{center}
{\em Instituto de Sistemas El\'etricos e Energia, Universidade
Federal de Itajub\'a\\Avenida BPS 1303, Pinheirinho, CEP
37.500-903, Itajub\'a, MG, Brazil
\\}e--mail:braga\_denis@yahoo.com.br
\end{center}

\begin{center}
{\bf Abstract}
\end{center}

In this paper we study the Lyapunov stability and the Hopf
bifurcation in a system coupling an hexagonal centrifugal governor
with a steam engine. Here are given sufficient conditions for the
stability of the equilibrium state and of the bifurcating periodic
orbit. These conditions are expressed in terms of the physical
parameters of the system, and hold for parameters outside a
variety of codimension two.

\vspace{0.1cm}



\newpage
\section{\bf Introduction}\label{intro}

The centrifugal governor is a device that automatically controls
the speed of an engine. The most important one, invented by James
Watt in 1788 --- Watt governor ---, is regarded as the starting
landmark for the theory of automatic control. The historical
relevance of this device as well as its importance for present day
theoretical and technological control developments --- going from
steam to diesel, gasoline engines and electronic governors ---
have been widely discussed by MacFarlane \cite{mac}, Denny
\cite{denny}, Fasol \cite{fasol} and Wellstead - Readman
\cite{control} among others.

The centrifugal governor design received several important
modifications as well as other types of governors were also
developed. From  MacFarlane \cite{mac}, p. 251, we quote:

\noindent ``Several important advances in automatic control
technology were made in the latter half of the 19th century. A key
modification to the flyball governor was the introduction of a
simple means of setting the desired running speed of the engine
being controlled by balancing the centrifugal force of the
flyballs against a spring, and using the preset spring tension to
set the running speed of the engine".

This paper is devoted to the study of the dynamic stability and
simplest bifurcations of the system coupling the hexagonal
centrifugal governor with a spring --- called  Hexagonal governor
--- and the steam engine. See Fig. \ref{hexagonal} for an
illustration. The system coupling the Hexagonal governor (resp.
Watt governor, with no spring and with  vanishing horizontal edges
of the hexagon) and the steam engine will be called simply the
Hexagonal Governor System (HGS) (resp. Watt Governor System
(WGS)). The stability analysis of the stationary states and of
small amplitude oscillations of this system will be pursued here.

The first mathematical analysis of the stability conditions in the
WGS was due to Maxwell \cite{max} and, in a user friendly style,
likely to be better understood by engineers, by Vyshnegradskii
\cite{vysh}. A simplified version of the WGS local stability based
on the work of Vyshnegradskii is presented by Pontryagin
\cite{pon}.

The oscillatory, small amplitude, behavior in the WGS has been
associated to a periodic orbit that appears from a Hopf
bifurcation. This was established by Hassard et al. in
\cite{has1}, Al-Humadi and Kazarinoff in \cite{humadi} and, in a
more general context, by the authors in \cite{smb1, smb2}.

In \cite{smb2}, restricting ourselves to Pontryagin's system of
differential equations for the WGS, we carried out a deeper
investigation of the stability of the equilibrium along the
critical Hopf bifurcations up to codimension 3, happening at a
unique point at which the bifurcation diagram was established. A
conclusion derived from the properties of the bifurcation diagram
implied the existence of parameters where the WGS has an
attracting periodic orbit coexisting  with an attracting
equilibrium.

The results of the present paper extend in a different direction
the analysis in \cite{smb1}, as described below.

In Section \ref{hexa} we introduce the differential equations that
model the HGS illustrated in Fig. \ref{hexagonal}. The stability
of the equilibrium point of this model is analyzed and a general
version of the stability condition is obtained and presented in
the terminology of Vyshnegradskii (Theorem \ref{teoestabilidade}
and Remark \ref{condVich}). The codimension 1 Hopf bifurcation for
the HGS differential equations is studied in Section \ref{hopf}.
An expression which determines the sign of the first Lyapunov
coefficient is obtained (Theorem \ref{teoremageral}). Sufficient
conditions for the stability of the bifurcating periodic orbit are
given. Two pertinent particular cases (no spring and vanishing
horizontal edge) are calculated and illustrated. See Theorem
\ref{teoremahopf}, Fig. \ref{sinalL1rho0} and Theorem
\ref{teoremahopf1}, Fig. \ref{sinalL1kappa0}.

Concluding comments are presented in Section \ref{conclusion}.

\section{The Hexagonal governor system}\label{hexa}

\newtheorem{teo}{Theorem}[section]
\newtheorem{lema}[teo]{Lemma}
\newtheorem{prop}[teo]{Proposition}
\newtheorem{cor}[teo]{Corollary}
\newtheorem{remark}[teo]{Remark}
\newtheorem{example}[teo]{Example}

\subsection{Hexagonal governor differential equations}\label{diffequat}

\begin{figure}[!h]
\centerline{
\includegraphics[width=12cm]{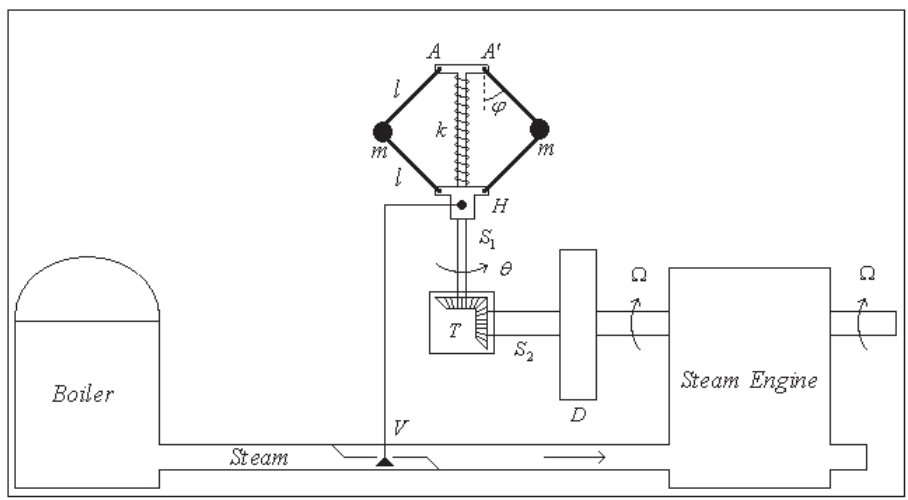}}
\caption{{\small Hexagonal centrifugal governor --- steam engine
system}.}

\label{hexagonal}
\end{figure}

The HGS studied in this paper is shown in Fig. \ref{hexagonal}.
There, $\varphi \in \left( 0,\frac{\pi}{2} \right)$ is the angle
of deviation of the arms of the governor from its vertical direction axis
$S_1$, $\Omega \in [0,\infty)$ is the angular velocity of the
rotation of the engine flywheel $D$, $\theta$ is the angular
velocity of the rotation of  $S_1$, $l$ is the length of the arms,
$m$ is the mass of each ball, $H$ is a  sleeve which supports the
arms and slides along $S_1$, $T$ is a set of transmission gears
and $V$ is the valve that determines the supply of steam to the
engine.

The HGS differential equations can be found as follows. For
simplicity, we neglect the mass of the sleeve and of the arms. There
are four forces acting on the balls at all times. They are the
tangential component of the gravity
\[
-m g \sin \varphi,
\]
where $g$ is the standard acceleration of gravity; the tangential
component of the centrifugal force
\[
m (L + l \sin \varphi ) \theta^2 \cos \varphi,
\]
$2L \geq 0$ is the distance $AA'$ in Fig. \ref{hexagonal}; the
tangential component of the restoring force due to the spring
\[
-2 k l (1 - \cos \varphi) \sin \varphi,
\]
$2l$ is the natural length of the spring and $k \geq 0$ is the
spring constant; and the force of friction
\[
- b l \dot {\varphi},
\]
$b > 0$ is the friction coefficient.

From Newton's Second Law of Motion, using the transmission
function $\theta = c \: \Omega$, where $c > 0$, one has
\begin{equation}
\ddot{\varphi} = c^2 \frac{L}{l} \: \Omega^2 \cos \varphi + \left(
\frac{2k}{m} + c^2 \Omega^2 \right)  \sin \varphi \cos \varphi -
\frac{2kl + mg}{ml} \: \sin \varphi -\frac{b}{m} \: \dot{\varphi}.
\label{Newton}
\end{equation}

The torque acting on the flywheel $D$ is
\begin{equation}
I \: \dot{\Omega} = \mu \: \cos \varphi - F, \label{torque}
\end{equation}
where $I$ is the moment of inertia of the flywheel, $F$ is an
equivalent torque of the load and $\mu > 0$ is a proportionality
constant to represent the torque due to the steam which decreases with the angle $\varphi$.
See \cite{pon}, p. 217, for more details.

From Eqs. (\ref{Newton}) and (\ref{torque}) the differential
equations of our model are given by
\begin{eqnarray}\label{sistema}
\frac{d \; \varphi}{d \tau} &=& \psi \nonumber\\
\frac{d \; \psi}{d \tau} &=& c^2 \frac{L}{l} \; \Omega^2 \cos
\varphi + \left( \frac{2k}{m} + c^2 \Omega^2 \right)  \sin \varphi
\cos \varphi -
\frac{2kl + mg}{ml} \; \sin \varphi -\frac{b}{m} \; \psi \nonumber \\
\frac{d \; \Omega}{d \tau} &=& \frac{1}{I} \; \left( \mu \cos
\varphi -F \right)
\end{eqnarray}
where $\tau$ is the time.

The standard Watt governor differential equations as presented in
Pontryagin \cite{pon}, p. 217, are obtained from (\ref{sistema})
by taking $L = 0$ and $k = 0$,
\begin{eqnarray}\label{standard}
\frac{d \; \varphi}{d \tau} &=& \psi \nonumber\\
\frac{d \; \psi}{d \tau} &=& c^2 \; \Omega^2 \; \sin \varphi \;
\cos \varphi -
\frac{g}{l} \; \sin \varphi - \frac{b}{m} \; \psi \\
\frac{d \; \Omega}{d \tau} &=& \frac{1}{I} \; \left( \mu \cos
\varphi -F \right) \nonumber
\end{eqnarray}

Performing the following changes in the coordinates, parameters and
time
\[
x = \varphi, \: y = \left( \frac{ml}{2kl+mg} \right)^{1/2} \psi,
\: z = c \left( \frac{ml}{2kl+mg} \right)^{1/2} \Omega,
\]
\begin{equation}
t = \left( \frac{2kl+mg}{ml} \right)^{1/2} \tau, \: \rho =
\frac{L}{l}, \: \kappa = \frac{2kl}{2kl+mg}, \label{novos}
\end{equation}
\[
\varepsilon = \frac{b}{m} \left( \frac{ml}{2kl+mg} \right)^{1/2},
\: \alpha = \frac{c \mu}{I} \left( \frac{ml}{2kl+mg} \right), \:
\beta = \frac{F}{\mu},
\]
where $\rho \geq 0$, $0 \leq \kappa < 1$, $\varepsilon > 0$,
$\alpha
> 0$ and $0 < \beta <1$, the differential equations
(\ref{sistema}) can be written as
\begin{eqnarray}\label{sistemafinal}
x' = \frac{d x}{d t} &=& y \nonumber\\
y' = \frac{d y}{d t} &=& \rho \; z^2 \; \cos x + (z^2 + \kappa) \;
\sin x \; \cos x -
\sin x - \varepsilon \; y \\
z' = \frac{d z}{d t} &=& \alpha \; (\cos x - \beta) \nonumber
\end{eqnarray}
or equivalently by
\begin{equation}
{\bf x}' = f({\bf x}, {\bf \zeta}),
\label{campo}
\end{equation}
where
\[
f({\bf x},{\bf \zeta}) = \left( y, \rho \; z^2 \; \cos x + (z^2 +
\kappa) \; \sin x \; \cos x - \sin x - \varepsilon \; y, \alpha \;
\left( \cos x - \beta \right) \right),
\]
\[
{\bf x} =(x,y,z) \in \left( 0, \frac{\pi}{2} \right) \times \R
\times [0,\infty)
\]
and
\[
{\bf \zeta} = (\beta, \alpha,\varepsilon, \rho, \kappa) \in \left(
0, 1\right) \times \left( 0, \infty \right) \times \left(0,\infty
\right) \times \left[ 0, \infty \right) \times \left[ 0, 1
\right).
\]

\subsection{Stability analysis at the equilibrium
point}\label{stability}

The HGS differential equations (\ref{sistemafinal}) have only one
admissible equilibrium point
\begin{equation}
P_0 = (x_0 , y_0 , z_0) = \left( \arccos \beta, 0, \frac{(1 -
\kappa \beta)^{1/2}(1 - \beta ^2)^{1/4}}{\beta ^{1/2} (\rho + (1 -
\beta^2)^{1/2})^{1/2}} \right). \label{P0}
\end{equation}

The Jacobian matrix of $f$ at $P_0$ has the form
\begin{equation}
Df \left( P_0 \right) = \left( \begin{array}{ccc}
0 & 1 & 0 \\
\\- \omega_0 ^2 & - \varepsilon & \xi \\
\\ -\alpha (1 - \beta^2)^{1/2}  & 0  & 0
\end{array} \right),
\label{jacobianP0}
\end{equation}
where
\begin{equation}
\omega_0 = \sqrt {\frac{(1 - \beta ^2)^{3/2} + \rho (1 - \kappa
\beta^3)}{\beta (\rho + (1 - \beta^2)^{1/2})}} \label{omega0}
\end{equation}
and
\[
\xi = 2 \beta^{1/2} (1 - \beta ^2)^{1/4} (1 - \kappa \beta)^{1/2}
(\rho + (1 - \beta ^2))^{1/2}.
\]

For the sake of completeness we state the following lemma whose
proof can be found in \cite{pon}, p. 58.

\begin{lema}
The polynomial $L(\lambda) = p_0 \lambda ^3 + p_1 \lambda ^2 + p_2
\lambda + p_3$, $p_0 > 0$, with real coefficients has all roots
with negative real parts if and only if the numbers $p_1 , p_2 ,
p_3$ are positive and the inequality $p_1 p_2 > p_0 p_3$ is
satisfied.

\label{routh}
\end{lema}

\begin{teo}
If
\begin{equation}
\varepsilon > \varepsilon_c = \frac{2 \alpha \beta^{3/2} (1 -
\beta ^2)^{3/4} (1 - \kappa \beta)^{1/2} (\rho + (1 - \beta
^2)^{1/2})^{3/2}}{(1 - \beta ^2)^{3/2} + \rho (1 - \kappa
\beta^3)}, \label{varepsiloncritico}
\end{equation}
then the HGS differential equations (\ref{sistemafinal}) have an
asymptotically stable equilibrium point at $P_0$. If
\[
0 < \varepsilon < \varepsilon_c
\]
then $P_0$ is unstable.

\label{teoestabilidade}
\end{teo}

\noindent {\bf Proof.} The characteristic polynomial of $Df \left(
P_0 \right)$ is given by $p(\lambda)$, where
\[
-p(\lambda)= \lambda^3 + p_1 \: \lambda^2 + p_2 \: \lambda + p_3,
\]
\[
p_1 = \varepsilon, \: \: p_2 = \frac{(1 - \beta^2)^{3/2} + \rho (1
- \kappa \beta^3)}{\beta (\rho + (1 - \beta^2)^{1/2})}
\]
and
\[
p_3 = \frac{2 \alpha \beta^{3/2} (1 - \beta ^2)^{3/4} (1 - \kappa
\beta)^{1/2} (\rho + (1 - \beta ^2)^{1/2})^{3/2}}{\beta (\rho + (1
- \beta^2)^{1/2} )}.
\]
The coefficients of $- p(\lambda)$ are positive. Thus a necessary
and sufficient condition for the asymptotic stability of the
equilibrium point $P_0$, as provided by the condition for one real
negative root and a pair of complex conjugate roots with negative
real part, is given by (\ref{varepsiloncritico}), according to
Lemma \ref{routh}.
\begin{flushright}
$\blacksquare$
\end{flushright}

\begin{remark}
In terms of the HGS physical parameters, condition
(\ref{varepsiloncritico}) is equivalent to
\begin{equation}
\frac{b \: I}{m} \: \eta > 1, \label{Vichneg}
\end{equation}
where
\begin{equation}
\eta = \left| \frac{d \Omega_0}{dF} \right| = \frac{(1 -
\beta^2)^{3/2} + \rho - \beta^3 \kappa \rho}{2 \beta^{3/2} (1 -
\beta^2)^{3/4} (1 - \kappa \beta)^{1/2} ((1 - \beta^2)^{1/2} +
\rho)^{3/2}} \label{naouniformidade}
\end{equation}
is the non-uniformity of the performance of the engine which
quantifies the change in the engine speed with respect to the load
(see \cite{pon}, p. 219, for more details). Eq.
(\ref{naouniformidade}) can be written in terms of the
original parameters of the HGS, but this expression is too long to be put in
print.

\label{condVich}
\end{remark}

The rules formulated by Vyshnegradskii to enhance the stability of
the system follow directly from (\ref{Vichneg}). In particular,
the interpretation of (\ref{Vichneg}) is that a sufficient amount
of damping ---$b$---  must be present relative to the other
physical parameters for the system to be stable at the desired
operating speed. Condition (\ref{Vichneg}) is equivalent to the
original condition given by Vyshnegradskii for the WGS (see
\cite{pon}, p. 219).

\section{Hopf bifurcation analysis}\label{hopf}

In this section we study the stability of $P_0$ under the
condition
\begin{equation}
\varepsilon = \varepsilon_c, \label{valorcritico}
\end{equation}
that is, on the Hopf hypersurface which is complementary to the range of validity of Theorem
\ref{teoestabilidade}.

\subsection{Generalities on Hopf bifurcations}

The study outlined below is based on the approach found in the
book of Kuznetsov \cite{kuznet}, pp 177-181.

Consider the differential equations
\begin{equation}
{\bf x}' = f ({\bf x}, {\bf \mu}), \label{diffequat}
\end{equation}
\noindent where ${\bf x} \in \R^3$ and ${\bf \mu} \in \R^m$ is a
vector of control parameters. Suppose (\ref{diffequat}) has an
equilibrium point ${\bf x} = {\bf x_0}$ at ${\bf \mu} = {\bf
\mu_0}$ and represent
\begin{equation}
F({\bf x}) = f ({\bf x}, {\bf \mu_0}) \label{Fhomo}
\end{equation}
as
\begin{equation}
F({\bf x}) = A{\bf x} + \frac{1}{2} \: B({\bf x},{\bf x}) +
\frac{1}{6} \: C({\bf x}, {\bf x}, {\bf x}) + O(|| {\bf x}
||^4){\nonumber}, \label{taylorexp}
\end{equation}
\noindent where $A = f_{\bf x}(0,{\bf \mu_0})$ and
\begin{equation}
B_i ({\bf x},{\bf y}) = \sum_{j,k=1}^3 \frac{\partial ^2
F_i(\xi)}{\partial \xi_j \: \partial \xi_k} \bigg|_{\xi=0} x_j \;
y_k, \label{Bap}
\end{equation}
\begin{equation}
C_i ({\bf x},{\bf y},{\bf z}) = \sum_{j,k,l=1}^3 \frac{\partial ^3
F_i(\xi)}{\partial \xi_j \: \partial \xi_k \: \partial \xi_l}
\bigg|_{\xi=0} x_j \; y_k \: z_l, \label{Cap}
\end{equation}
\noindent for $i = 1, 2, 3$. Here the variable ${\bf x}-{\bf x_0}$
is also denoted by ${\bf x}$.

Suppose $({\bf x_0}, {\bf \mu_0})$ is an equilibrium point of
(\ref{diffequat}) where the Jacobian matrix $A$ has a pair of
purely imaginary eigenvalues $\lambda_{2,3} = \pm i \omega_0$,
$\omega_0 > 0$, and no other critical (i.e., on the imaginary
axis) eigenvalues.

The two dimensional center manifold can be parametrized by $w \in
\R^2 = \C$,  by means of ${\bf x} = H (w,\bar w )$, which is
written as
\[
H(w,{\bar w}) = w q + {\bar w}{\bar q} + \sum_{2 \leq j+k \leq 3}
\frac{1}{j!k!} \: h_{jk}w^j{\bar w}^k + O(|w|^4),
\]
with $h_{jk} \in \C ^3$, $h_{jk}={\bar h}_{kj}$.

Substituting these expressions into (\ref{diffequat}) and
(\ref{taylorexp}) one has
\begin{equation}
H_w (w,\bar w )w' + H_{\bar w}  (w,\bar w ){\bar w}'  = F(H(w
,\bar w )).
\label{homologicalp}
\end{equation}

Let $p, q \in \C ^3$ be vectors such that
\begin{equation}
A q = i \omega_0 \: q,\:\: A^{\top} p = -i \omega_0 \: p, \:\:
\langle p,q \rangle = \sum_{i=1}^3 \bar{p}_i \: q_i \:\: = 1.
\label{normalization}
\end{equation}

The complex vectors $h_{ij}$ are to be determined so that equation
(\ref{homologicalp}) writes as follows
\begin{equation}
w'= i \omega_0 w + \frac{1}{2} \: G_{21} w |w|^2 + O(|w|^4),
\label{omegalinha}
\end{equation}
with  $G_{21} \in \C $.

Solving the linear system obtained by expanding
(\ref{homologicalp}), the coefficients of the quadratic terms of
(\ref{Fhomo}) lead to
\begin{equation}
h_{11}=-A^{-1}B(q,{\bar q}) \label{h11},
\end{equation}
\begin{equation}
h_{20}=(2i\omega_0 I_3 - A)^{-1}B(q,q),\label{h20}
\end{equation}
where $I_3$ is the unit $3 \times 3$ matrix.

The coefficients of the cubic terms are also uniquely calculated,
except for the term $w^2 {\bar w}$, whose coefficient satisfies a
singular system for $h_{21}$
\begin{equation}
(i \omega_0 I_3 -A)h_{21}=C(q,q,{\bar q})+B({\bar q},h_{20}) + 2
B(q,h_{11})-G_{21}q, \label{h21m}
\end{equation}
which has a solution if and only if
\[
\langle p, C(q,q,\bar q) + B(\bar q, h_{20}) + 2 B(q,h_{11})
-G_{21} q \rangle = 0.
\]
Therefore
\begin{equation} \label{G21}
G_{21}= \langle p, C(q,q,\bar q) + B(\bar q, (2i \omega_0
I_3-A)^{-1} B(q,q)) - 2 B(q,A^{-1} B(q,\bar q)) \rangle.
\end{equation}
The {\it first Lyapunov coefficient} $l_1$ is defined
by
\begin{equation}
l_1 =  \frac{1}{2 \; \omega_0} \: {\rm Re} \; G_{21}.
\label{defcoef}
\end{equation}

From (\ref{omegalinha}) its sign decides the stability, when
negative, or instability, when positive, of the equilibrium.

A {\it Hopf point} $({\bf x_0}, {\bf \mu_0})$ is an equilibrium
point of (\ref{diffequat}) where the Jacobian matrix $A$ has a
pair of purely imaginary eigenvalues $\lambda_{2,3} = \pm i
\omega_0$, $\omega_0 > 0$, and no other critical eigenvalues. At a
Hopf point, a two dimensional center manifold  is well-defined,
which is invariant under the flow generated by (\ref{diffequat})
and can be smoothly continued to nearby parameter values.

A Hopf point is called {\it transversal} if the curves of complex
eigenvalues cross the imaginary axis with non-zero derivative.

In a neighborhood of a transversal Hopf point with $l_1 \neq 0$
the dynamic behavior of the system (\ref{diffequat}), reduced to
the family of parameter-dependent continuations of the center
manifold, is orbitally topologically equivalent to the complex
normal form
\begin{equation}\label{nf}
w' = (\gamma + i \omega) w + l_1 w |w|^2 ,
\end{equation}
$w \in \C $, $\gamma$, $\omega$ and $l_1$ are smooth continuations
of $0$, $\omega_0$ and the first Lyapunov coefficient at the Hopf
point \cite{kuznet}. When $l_1 < 0$ ($l_1
> 0$) a family of stable (unstable) periodic orbits can be found
on this family of center manifolds, shrinking to the equilibrium
point at the Hopf point.

\subsection{Hopf bifurcation in the HGS}

From (\ref{campo}) write the Taylor expansion (\ref{taylorexp}) of
$f({\bf x})$. Define
\begin{equation}
\omega_1 = \sqrt {\frac{1 - \beta ^2}{\beta}} \label{omega1}
\end{equation}
and
\begin{equation}
\sigma = \sqrt {\frac{1 - \kappa \beta}{\rho + \omega_1
\beta^{1/2}}}. \label{sigma}
\end{equation}

Thus, with $\omega_0$ given in Eq. (\ref{omega0}),
\begin{equation}
A = \left( \begin{array}{ccc}
0 & 1 & 0 \\
\\- \omega_0 ^2 & - \varepsilon_c & \displaystyle \frac {\varepsilon_c \; \omega_0 ^2}{\alpha \beta ^{1/2} \omega_1} \\
\\ -\alpha \beta ^{1/2} \omega_1    & 0  & 0
\end{array} \right),
\label{partelinear}
\end{equation}
and, with the notation in (\ref{taylorexp}) we have
\begin{equation}
F({\bf x})\, - \, A{\bf x} = \left( 0,F_2({\bf x})+ O(||x||^4),
F_3({\bf x}) + O(||x||^4) \right), \label{partenaolinear}
\end{equation}
where
\begin{eqnarray*}
F_2({\bf x}) = - \frac{3}{2} \beta^{1/2} \omega_1 (1 - \rho \sigma
^2) x^2 + \frac{2 \sigma (\beta^{1/2} \omega_1)^{1/2} (2 \beta^2
-1 - \beta^{1/2} \rho \omega_1)}{\beta^{1/2}} \; x z + \\ \beta
(\rho + \beta^{1/2} \omega_1) z^2 + \frac{1 + (3 - 7 \beta^2)(1 -
\rho \sigma^2)}{6 \beta} \; x^3 - \\ \beta^{1/2} \sigma
(\beta^{1/2} \omega_1)^{1/2} (\rho + 4 \beta^{1/2} \omega_1) \;
x^2 z + (2 \beta^2 -1 - \beta^{1/2} \rho \omega_1) \; x z^2,
\end{eqnarray*}
and
\begin{eqnarray*}
F_3({\bf x}) = -\frac{1}{2} \; \alpha \; \beta \; x^2 +
\frac{1}{6} \; \alpha \; \beta^{1/2} \omega_1 \; x^3.
\end{eqnarray*}

From (\ref{partelinear}) the eigenvalues of $A$ are
\begin{equation}
\lambda_1 = -\varepsilon_c , \: \: \lambda_2 = i \: \omega_0, \:
\: \lambda_3 = -i \: \omega_0. \label{autovalores}
\end{equation}
\noindent The eigenvectors $q$ and $p$ satisfying
(\ref{normalization}) are respectively
\begin{equation}
q = \left( -i, \omega_0, \frac{\alpha \beta^{1/2}
\omega_1}{\omega_0} \right) \label{q}
\end{equation}
\noindent and
\begin{equation}
p = \left( -\frac{i}{2}, \frac{\omega_0 - i \varepsilon_c}{2 (
\omega_0 ^2 + \varepsilon_c ^2 )}, \frac{\varepsilon_c \omega_0
(\varepsilon_c + i \omega_0)}{ 2 \alpha \beta^{1/2} \omega_1
(\omega_0 ^2 + \varepsilon_c ^2) } \right). \label{p}
\end{equation}

The main result of this section can be formulated now.

\begin{teo}
Consider the family of differential equations
(\ref{sistemafinal}). The first Lyapunov coefficient at the point
(\ref{P0}) for parameter values satisfying (\ref{valorcritico}) is
given by
\begin{equation}
l_1 (\beta, \alpha, \rho, \kappa) = - \frac{R(\beta, \alpha,
 \rho, \kappa)}{4 \beta \varepsilon_c \omega_0 ^4
\omega_1 ^2 (\varepsilon_c ^4 + 5 \varepsilon_c ^2 \omega_0 ^2 + 4
\omega_0 ^4)} ,\label{coeficiente1}
\end{equation}
\noindent where
{\small
\begin{eqnarray*}
R(\beta, \alpha, \rho, \kappa) = \varepsilon_c ^2 \omega_0 ^4
\omega_1 ^2 (3 \rho \sigma^2 -4)(\varepsilon_c ^2 + 4 \omega_0 ^2)
+ \beta \varepsilon_c ^2 \omega_0 ^6 \omega_1 ^2 (\varepsilon_c ^2
+ 4 \omega_0 ^2) + \\ 8 \alpha \beta ^{11/4} \varepsilon_c \sigma
\omega_0 ^4 \omega_1 ^{9/2} (\varepsilon_c ^2 + 4 \omega_0 ^2)+ 8
\alpha^2 \beta ^{3/2} \varepsilon_c^2 \sigma^2 \omega_0 ^2
\omega_1 ^5 - \\ 4 \alpha ^3 \beta ^{13/4} \varepsilon_c \rho
\sigma \omega_1 ^{11/2} (\varepsilon_c ^2 + 10 \omega_0 ^2) + 8
\alpha^3 \beta ^{21/4} \varepsilon_c \rho \sigma \omega_1 ^{11/2}
(\varepsilon_c ^2 + 10 \omega_0 ^2) + \\ 4 \alpha^4 \beta^5 \rho^2
\omega_1 ^6 (\varepsilon_c ^2 + 8 \omega_0 ^2) + 8 \alpha^3
\beta^{23/4} \varepsilon_c \sigma \omega_1 ^{13/2} (\varepsilon_c
^2 + 10 \omega_0 ^2) - \\ 4 \alpha^3 \beta^{17/4} \varepsilon_c
\rho \sigma \omega_1 ^{15/2} (\varepsilon_c ^2 + 10 \omega_0 ^2) +
4 \alpha^4 \beta^6 \omega_1 ^8 (\varepsilon_c ^2 + 8 \omega_0 ^2)
- \\ 2 \alpha \beta^{7/4} \varepsilon_c \sigma \omega_0 ^2
\omega_1 ^{5/2} (14 \omega_0 ^4 + 3 \varepsilon_c ^2 \omega_1 ^2
(\rho \sigma ^2 -1) + 30 \omega_0 ^2 \omega_1 ^2 (\rho \sigma ^2
-1)) + \\ 8 \beta^{11/2} (4 \alpha^2 \varepsilon_c ^2 \sigma^2
\omega_0 ^2 \omega_1 ^5 +
\alpha^4 \rho \omega_1 ^7 (\varepsilon_c ^2 + 8 \omega_0 ^2)) - \\
4 \alpha \beta^{15/4} \varepsilon_c \sigma \omega_1 ^{5/2} (-14
\omega_0 ^6 - 30 \omega_0 ^4 \omega_1 ^2 (\rho \sigma ^2 -1) + \\
\alpha ^2 \varepsilon_c ^2 \omega_1 ^4 (1 + \rho ^2) + \omega_0 ^2
\omega_1 ^2 ( 3 \varepsilon_c ^2 ( 1 - \rho \sigma ^2) + 10 \alpha
^2 \omega_ 1 ^2 (1 +\rho ^2))) + \\ 2 \alpha^2 \beta ^{5/2}
\varepsilon_c ^2 \rho \omega_0 ^2 \omega_1 ^5 (\varepsilon_c ^2 +
4 (\omega_0 ^2 + \rho \sigma ^2 \omega_1 ^2)) + 2 \alpha \beta
^{9/4} \varepsilon_c \rho \sigma \omega_0 ^2 \omega_1 ^{7/2} \\
(\varepsilon_c ^2 (\omega_0 ^2 - 3 \omega_1 ^2 (\rho \sigma ^2
-1)) - 10 \omega_0 ^2 (\omega_0 ^2 + 3 \omega_1 ^2 (\rho \sigma ^2
-1))) - \\ 2 \alpha ^2 \beta ^4 \omega_0 ^2 \omega_1 ^4 (2
\varepsilon_c ^4 - 8 \omega_0 ^2 (\omega_0 ^2 + 3 \omega_1 ^2
(\rho \sigma ^2 -1))+ \varepsilon_c ^2 (9 \omega_0 ^2 + \omega_1
^2 (3 + 13 \rho \sigma ^2))) + \\ \beta^2 \varepsilon_c ^2
\omega_0 ^2 (\varepsilon_c ^2 (\omega_0 ^4 - 4 \omega_0 ^2
\omega_1 ^2 (\rho \sigma ^2 -1) + 2 \alpha^2 \omega_1 ^4 ) + 2 (4
\omega_0 ^6 + \omega_0 ^4 \omega_1 ^2 (\rho \sigma ^2 -1) +\\
\omega_0 ^2 \omega_1 ^4 (4 \alpha^2 + 9 (\rho \sigma ^2 -1)^2) + 8
\alpha ^2 \rho \sigma^2 \omega_1 ^6)) + \\ 2 \alpha ^2 \beta
^{7/2} \omega_0 ^2 \omega_1 ^3 (8 \rho \omega_0 ^2 (\omega_0 ^2 +
3 \omega_1 ^2 (\rho \sigma ^2 -1)) + \\ \varepsilon_c ^2 (-16
\sigma^2 \omega_1 ^2 + 3 \rho ^2 \sigma ^2 \omega_1 ^2 - \rho
(\omega_0 ^2 + 3 \omega_1 ^2))).
\end{eqnarray*}
}
\label{teoremageral}
\end{teo}

\noindent{\bf Proof.} The proof depends on preliminary
calculations presented below. From (\ref{taylorexp}), (\ref{Bap}),
(\ref{Cap}) and (\ref{partenaolinear}) one has
\begin{equation}
B({\bf x},{\bf y}) = \left( 0, B_2({\bf x},{\bf y}), - \alpha
\beta \: x_1 \: y_1 \right), \label{B1}
\end{equation}
\noindent where
\begin{eqnarray*}
B_2({\bf x},{\bf y})= -3 \: \beta^{1/2} \omega_1 (1 - \rho
\sigma^2) \: x_1 \: y_1 + 2 \beta (\rho + \omega_1 \beta ^{1/2})
\: x_3 \: y_3 + \\ \frac{2 \sigma (\beta^{1/2} \omega_1)^{1/2}((2
\beta^2 -1)- \rho \omega_1 \beta^{1/2})}{\beta^{1/2}} \: (x_1 \:
y_3 + x_3 \: y_1),
\end{eqnarray*}
\begin{equation}
C({\bf x},{\bf y}, {\bf z}) = \left( 0, C_2({\bf x},{\bf y},{\bf
z}), \alpha \beta^{1/2} \omega_1 \: x_1 \: y_1 \: z_1 \right),
\label{C1}
\end{equation}
where
\begin{eqnarray*}
C_2({\bf x},{\bf y},{\bf z})= \frac{1 + (1- \rho \sigma^2)(3-7
\beta^2)}{\beta} \: x_1 \: y_1 \: z_1 + 2 (2 \beta^2 -1 - \rho
\beta^{1/2} \omega_1) \\ (x_1 y_3 z_3 + x_3 y_1 z_3 + x_3 y_3 z_1)
- 2 \beta^{1/2} \sigma (\beta^{1/2} \omega_1)^{1/2}(\rho + 4
\beta^{1/2} \omega_1) \\ (x_1 y_1 z_3 + x_1 y_3 z_1 + x_3 y_1
z_1).
\end{eqnarray*}

Referring to the notation in (\ref{B1}), (\ref{C1}) and (\ref{q})
one has
\begin{equation}
B(q,q) = \left( 0, B_2(q,q), \alpha \beta \right), \label{Bqq}
\end{equation}
where
\begin{eqnarray*}
B_2(q,q) = \frac{\beta \omega_1 ^2}{\omega_0 ^2 (\beta^{1/2}
\omega_1)^{3/2}} \Bigg[ 2 \alpha^2 \beta (\beta^{1/2}
\omega_1)^{3/2} (\rho + \beta^{1/2} \omega_1) + \\ 3 (1 - \rho
\sigma^2) \omega_0 ^2 (\beta^{1/2} \omega_1)^{1/2} + i 4 \alpha
\sigma \omega_0 \omega_1 (1 - 2 \beta^2 + \beta^{1/2} \rho
\omega_1) \Bigg],
\end{eqnarray*}
\begin{equation}
B(q,\bar{q}) = \left( 0, B_2(q,\bar q),  - \alpha \beta \right),
\label{Bqqbar}
\end{equation}
where
\[
B_2(q,\bar q) = \frac{ \beta^{1/2} \omega_1 (3 \omega_0 ^2 (\rho
\sigma^2 -1) + 2 \alpha^2 \beta^{3/2} \omega_1 (\rho + \beta^{1/2}
\omega_1))}{\omega_0 ^2},
\]
\begin{equation}
C(q,q,\bar{q}) = \left( 0,  C_2(q,q,\bar q), -i \alpha \beta^{1/2}
\omega_1 \right), \label{Cqqqbar}
\end{equation}
where
\begin{eqnarray*}
C_2(q,q,\bar q)= \frac{-i}{\beta \omega_0 ^2} \Bigg[ \omega_0 ^2
(4 - 3 \rho \sigma^2 + 7 \beta^2 (\rho \sigma^2 -1))+ 2 \alpha^2
\beta^2 \omega_1 ^2 \\ (2 \beta^2 -1 -\beta^{1/2} \rho \omega_1) -
i 2 \alpha \beta^2 \sigma \omega_0 \omega_1 (\beta^{1/2}
\omega_1)^{1/2} (\rho + 4 \beta^{1/2} \omega_1) \Bigg].
\end{eqnarray*}
The first Lyapunov coefficient is given by (\ref{defcoef}). From
(\ref{p}) and (\ref{Cqqqbar}) one has
\begin{eqnarray}\label{partereal1}
{\rm Re} \langle p, C(q,q,\bar{q}) \rangle = \frac{-1}{2 \beta
\omega_0 ^2 (\varepsilon_c ^2 + \omega_0 ^2)} \Bigg[ 2 \alpha
\beta^{9/4} \sigma \omega_0 ^2 \omega_1 ^{3/2} (\rho + 4
\beta^{1/2} \omega_1) + \nonumber \\ + \varepsilon_c \Bigg(
\omega_0 ^2 (3 \rho \sigma^2 -4 + 7 \beta^2 (1 - \rho \sigma^2)) +
\beta \omega_0 ^4 + \\ 2 \alpha^2 \beta^2 \omega_1 ^2 (1 - 2
\beta^2 + \beta^{1/2} \rho \omega_1) \Bigg) \Bigg] \nonumber.
\end{eqnarray}
\noindent From (\ref{p}), (\ref{B1}), (\ref{q}) and (\ref{h11})
one has
\begin{eqnarray}\label{partereal2}
{\rm Re} \langle p, 2 B(q,h_{11}) \rangle =
\frac{-\beta^{3/4}}{\varepsilon_c \omega_0 ^4 \omega_1 ^2
(\varepsilon_c ^2 + \omega_0 ^2)} \Bigg[ 4 \alpha ^3 \varepsilon_c
\rho \sigma \omega_1 ^{11/2} (2 \beta^{7/2} - \beta^{3/2})+
\nonumber \\ 4 \alpha ^4 \beta^{13/4} \rho^2 \omega_1 ^6 + 8
\alpha ^3 \beta ^4 \varepsilon_c \sigma \omega_1 ^{13/2} + 8
\alpha ^4 \beta^{15/4} \rho \omega_1 ^7 - 4 \alpha ^3 \beta^{5/2}
\varepsilon_c \rho \sigma \omega_1^{15/2} + \nonumber \\ 4 \alpha
^4 \beta^{17/4} \omega_1 ^8 + \beta^{1/4} \varepsilon_c ^2
\omega_0 ^4 (\omega_0 ^2 + 3 \omega_1 ^2 (\rho \sigma ^2 -1))+
\nonumber \\ 2 \alpha ^2 \beta^{7/4} \rho \omega_0 ^2 \omega_1 ^3
(\omega_0 ^2 + 3 \omega_1 ^2 (\rho \sigma ^2 -1)) + \nonumber \\ 2
\alpha ^2 \beta^{9/4} \omega_0 ^2 \omega_1 ^4 (\omega_0 ^2 + 3
\omega_1 ^2 (\rho \sigma ^2 -1)) - \nonumber \\ 2 \alpha
\varepsilon_c \sigma \omega_0 ^2 \omega_1 ^{5/2} (2 \omega_0 ^2 +
3 \omega_1 ^2 (\rho \sigma ^2 -1)) - \\ 2 \alpha \beta ^{1/2}
\varepsilon_c \rho \sigma \omega_0 ^2 \omega_1 ^{7/2} (2 \omega_0
^2 + 3 \omega_1 ^2 (\rho \sigma ^2 -1)) - \nonumber \\ 4 \alpha
\beta ^2 \varepsilon_c \sigma \omega_1 ^{5/2} (\alpha ^2 \omega_1
^4 (1 + \rho ^2) -3 \omega_0 ^2 \omega_1 ^2 (\rho \sigma ^2 -1) -
2 \omega_0 ^4) \Bigg] \nonumber.
\end{eqnarray}
\noindent From (\ref{p}), (\ref{B1}), (\ref{q}) and (\ref{h20})
one has
\begin{equation}
{\rm Re} \langle p, B(\bar{q},h_{20}) \rangle =
\frac{\vartheta(\beta, \alpha, \rho, \kappa)}{2 \omega_0 ^4
\omega_1 ^2 (\varepsilon_c ^4 + 5 \varepsilon_c ^2 \omega_0 ^2 + 4
\omega_0 ^4)}, \label{partereal3}
\end{equation}
where
\begin{eqnarray*}
\vartheta(\beta, \alpha, \rho, \kappa)= - 8 \alpha^2 \beta^{1/2}
\varepsilon_c \sigma ^2 \omega_0 ^2 \omega_1 ^5 + 4 \alpha ^3
\sigma \rho \omega_1 ^{11/2} (\varepsilon_c ^2 - 2 \omega_0 ^2) (2
\beta ^{17/4}- \beta ^{9/4}) + \\ 4 \alpha^4 \beta^4 \varepsilon_c
\sigma^2 \omega_1 ^6 + 8 \alpha ^3 \beta^{19/4} \sigma \omega_1
^{3/2} (\varepsilon_c ^2 - 2 \omega_0 ^2) - \\ 8 \alpha^2
\beta^{3/2} \varepsilon_c \sigma^2 \rho^2 \omega_0 ^2 \omega_1 ^7
- 4 \alpha^3 \beta^{13/4} \sigma \rho \omega_1 ^{15/2}
(\varepsilon_c ^2 - 2 \omega_0 ^2) + \\ 4 \alpha ^4 \beta ^5
\varepsilon_c \omega_1 ^8 + 8 \alpha ^2 \beta^{9/2} \varepsilon_c
\omega_1 ^5 (\alpha ^2 \rho \omega_1 ^2 - 4 \sigma^2 \omega_0 ^2)
+ \\ 2 \alpha ^2 \beta ^3 \varepsilon_c \omega_0 ^2 \omega_1 ^4 (3
\omega_0 ^2 + \omega_1 ^2 (19 \rho \sigma ^2 -3 )) + \\ 2 \alpha
^2 \beta ^{5/2} \varepsilon_c \omega_0 ^2 \omega_1 ^3 (16 \sigma
^2 \omega_1 ^2 + 3 \rho ^2 \sigma ^2 \omega_1 ^2 + 3 \rho
(\omega_0 ^2 - \omega_1 ^2))- \\ 2 \alpha \beta ^{3/4} \sigma
\omega_0 ^2 \omega_1 ^{5/2} (2 \omega_0 ^2 (\omega_0 ^2 - 3
\omega_1 ^2 (\rho \sigma ^2 -1))+ \varepsilon_c ^2 (4 \omega_0 ^2
+ 3 \omega_1 ^2 (\rho \sigma ^2 -1)))- \\ 2 \alpha \beta ^{5/4}
\rho \sigma \omega_0 ^2 \omega_1 ^{7/2} (2 \omega_0 ^2 (\omega_0
^2 - 3 \omega_1 ^2 (\rho \sigma ^2 -1))+ \varepsilon_c ^2 (4
\omega_0 ^2 + 3 \omega_1 ^2 (\rho \sigma ^2 -1)))- \\ 4 \alpha
\beta ^{11/4} \sigma \omega_1 ^{5/2} (\varepsilon_c ^2 (-4
\omega_0 ^4  - 3 \omega_0 ^2 \omega_1 ^2 (\rho \sigma ^2 -1) +
\alpha ^2 \omega_1 ^4 (1 + \rho ^2)) - \\ 2 \omega_0 ^2 (\omega_0
^4 - 3 \omega_0 ^2 \omega_1 ^2 (\rho \sigma ^2 -1) + \alpha^2
\omega_1 ^4 (1 + \rho ^2))) + \\ \beta \varepsilon_c \omega_0 ^2
(\varepsilon_c ^2 \omega_0 ^2 \omega_1 ^2 (\omega_0 ^2 + 3 (\rho
\sigma ^2 -1)) - \\ 2 \omega_1 ^2 (3 \omega_0 ^4 (\rho \sigma^2
-1) + 9 \omega_0 ^2 \omega_1 ^2 (\rho \sigma ~2 -1) + 8 \alpha^2
\rho \sigma ^2 \omega_1 ^4)).
\end{eqnarray*}
\noindent Substituting (\ref{partereal1}), (\ref{partereal2}) and
(\ref{partereal3}) into (\ref{G21}) and (\ref{defcoef}), the
theorem is proved.
\begin{flushright}
$\blacksquare$
\end{flushright}

\begin{remark}\label{denominador}
The denominator of the first Lyapunov coefficient given by Eq.
(\ref{coeficiente1}) is positive. Thus the sign of $l_1$ is
determined by the sign of the function $R$, the numerator of
$l_1$.

The expression for $l_1$ depends only the parameters $\alpha,
\beta, \rho$ and $ \kappa$, although in the expression in
(\ref{coeficiente1}) appear also $\omega_0, \omega_1, \sigma $ and
$\varepsilon_c$.  This is due to the fact that these last
parameters are functions of the previous ones as shown in
(\ref{omega0}), (\ref{omega1}) and (\ref{sigma}).
\end{remark}

\begin{prop}
Consider the family of differential equations (\ref{sistemafinal})
regarded as dependent on the parameter $\varepsilon$. The real
part, $\gamma$, of the pair of complex eigenvalues verifies
\begin{equation}
\gamma'(\varepsilon_c) = - \frac{\omega_0 ^2}{2 (\omega_0 ^2 +
\varepsilon_c ^2)} < 0. \label{transversal}
\end{equation}
Therefore, the transversality condition holds at the Hopf point.
\label{lematransv}
\end{prop}

\noindent{\bf Proof.} Let $\lambda (\varepsilon) = \lambda_{2,3}
(\varepsilon) = \gamma (\varepsilon) \pm i \omega (\varepsilon)$
be eigenvalues of $A(\varepsilon)$ such that $\gamma
(\varepsilon_c) = 0$ and $\omega (\varepsilon_c) = \omega_0$,
according to (\ref{autovalores}). Taking the inner product of $p$
with the derivative of $A(\varepsilon) q(\varepsilon) = \lambda
(\varepsilon) q(\varepsilon)$ at $\varepsilon = \varepsilon_c$ one
has
\[
\left \langle p , \frac{d A}{d \varepsilon} \Bigg |_ {\varepsilon
= \varepsilon_c} \: q \right \rangle = \gamma'(\varepsilon_c) \pm
\omega'(\varepsilon_c).
\]
Thus the transversality condition is given by
\begin{equation}
\gamma'(\varepsilon_c) = {\rm Re} \: \left \langle p, \frac{d A}{d
\varepsilon} \Bigg |_ {\varepsilon = \varepsilon_c} \: q \right
\rangle. \label{deftransv}
\end{equation}
As
\[
\frac {d A}{d \varepsilon} \Bigg |_ {\varepsilon = \varepsilon_c}
\: q = \left( 0, - \omega_0, 0 \right),
\]
the proposition follows from (\ref{p}) and a simple calculation.
\begin{flushright}
$\blacksquare$
\end{flushright}

The full expression of $l_1$ in terms of the parameters $\alpha,
\beta, \rho, \kappa$ seems too long to be of use in qualitative
arguments. Two special cases are considered below for the sake of
illustration.

\subsubsection{The case $\rho = 0$}\label{rho=0}

\begin{cor}
Consider the case where $\rho = 0$. Then the equilibrium point
$P_0$ in (\ref{P0}) is given by
\begin{equation}
P_0 = (x_0 , y_0 , z_0) = \left( \arccos \beta, 0, \left (
\frac{1}{\beta } - \kappa  \right)^{1/2} \right), \label{P0rho=0}
\end{equation}
the Hopf hypersurface (\ref{valorcritico}) is given by
\begin{equation}
\varepsilon_c = \varepsilon_c (\beta, \alpha, \kappa) = 2 \;
\alpha \; \beta^{3/2} \; (1 - \kappa \beta)^{1/2}
\label{valorcriticorho=0}
\end{equation}
and the numerator of $l_1$ in (\ref{coeficiente1}) is given by
\begin{eqnarray}\label{l1rho=0}
G_1(\beta, \alpha, \kappa) = -3 + 5 \kappa \beta - (\alpha^2 -5)
\beta^2 + \kappa (\alpha ^2 -7) \beta^3 - \nonumber \\  2 \alpha^2
\kappa^2 \beta^4 - (\alpha^4 - 2 \alpha^2 \kappa^2) \beta ^6  +
\alpha ^4 \kappa \beta^7 .
\end{eqnarray}
If $G_1$ is different from zero then the family of HGS
differential equations (\ref{sistemafinal}) has a transversal Hopf
point at $P_0$ for $\varepsilon_c = 2 \; \alpha \; \beta^{3/2} \;
(1 - \kappa \beta)^{1/2}$.

\label{teoremarho=0}
\end{cor}

\noindent {\bf Proof.} The proof is immediate by substituting
$\rho = 0$ into Eqs. (\ref{P0}), (\ref{valorcritico}) and
(\ref{coeficiente1}). A sufficient condition for being a Hopf
point is that the first Lyapunov coefficient $l_1 \neq 0$, since
the transversality condition is satisfied by Proposition
\ref{lematransv}. But from (\ref{l1rho=0}) it is equivalent to
$G_1 \neq 0$.
\begin{flushright}
$\blacksquare$
\end{flushright}

\begin{remark}
The expression (\ref{P0rho=0}) shows that the ``running speed" of
the system depends monotonically decreasing on $\kappa$, which is
monotonically increasing on $k$, according to (\ref{novos}). This
corroborates analytically the quotation of MacFarlane in the
Introduction.

\label{velocidade}
\end{remark}

\begin{figure}[!h]
\centerline{
\includegraphics[width=10cm]{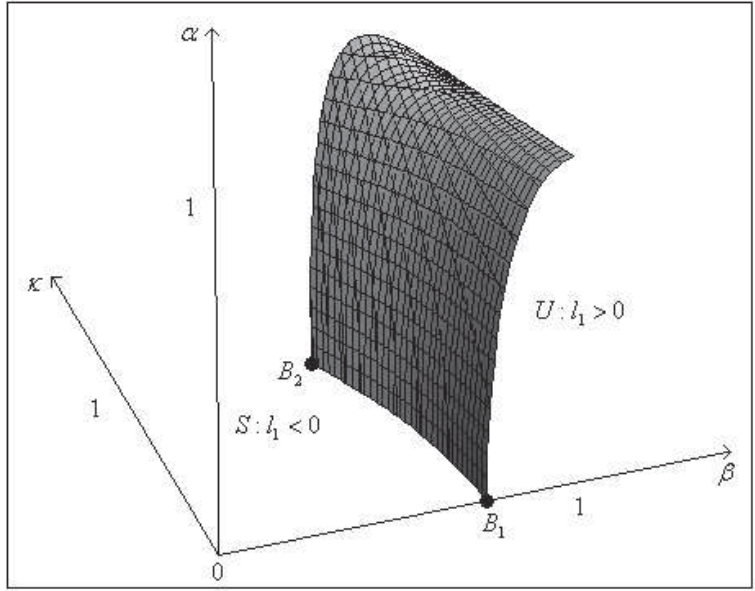}}

\caption{{\small Signs of the first Lyapunov coefficient for $\rho
= 0$}.}

\label{sinalL1rho0}
\end{figure}

Equation (\ref{l1rho=0}) gives a simple expression to determine
the sign of the first Lyapunov coefficient (\ref{coeficiente1})
for the case $\rho = 0$. The graph $G_1(\beta, \alpha,\kappa) = 0$
is illustrated in Fig. \ref{sinalL1rho0}, where the signs of the
first Lyapunov coefficient are also represented. The surface $l_1
= 0$ divides the hypersurface of critical parameters
$\varepsilon_c = 2 \; \alpha \; \beta^{3/2} \; (1 - \kappa
\beta)^{1/2}$ into two connected components denoted by $S$ and $U$
where $l_1 < 0$ and $l_1 > 0$ respectively. In  Fig.
\ref{sinalL1rho0}, the $\beta$ coordinates at the reference points
$B_1$ and $B_2$ are $0.7746$ and $0.5272$, respectively.

The following theorem summarizes the results in this subsection.

\begin{teo}
Consider the case where $\rho = 0$. If $(\beta, \alpha, \kappa)
\in S \cup U$ then the family of differential equations
(\ref{sistemafinal}) has a transversal Hopf point at $P_0$ for
$\varepsilon = \varepsilon_c$. If $(\beta, \alpha, \kappa) \in S$
then the Hopf point at $P_0$ for $\varepsilon = \varepsilon_c$ is
asymptotically stable and for each $\varepsilon < \varepsilon_c$,
but close to $\varepsilon_c$, there exists a stable periodic orbit
near the unstable equilibrium point $P_0$. If $(\beta, \alpha,
\kappa) \in U$ then the Hopf point at $P_0$ for $\varepsilon =
\varepsilon_c$ is unstable and for each $\varepsilon
> \varepsilon_c$, but close to $\varepsilon_c$, there exists an
unstable periodic orbit near the asymptotically stable equilibrium
point $P_0$. See Fig \ref{sinalL1rho0}.

\label{teoremahopf}
\end{teo}

\subsubsection{The case $\kappa = 0$}\label{kappa=0}

\begin{cor}
Consider the case where $\kappa = 0$. Then the equilibrium point
$P_0$ in (\ref{P0}) is given by
\begin{equation}
P_0 = (x_0 , y_0 , z_0) = \left( \arccos \beta, 0,
\frac{(1-\beta^2)^{1/4}}{\beta^{1/2} (\rho + (1 -
\beta^2)^{1/2})^{1/2}} \right), \label{P0kappa=0}
\end{equation}
the Hopf hypersurface (\ref{valorcritico}) is given by
\begin{equation}
\varepsilon_c = \varepsilon_c (\beta, \alpha, \rho) = \frac{2
\alpha \beta^{3/2} \; (1 - \beta^2)^{3/4} (\rho + (1 -
\beta^2)^{1/2})^{3/2}}{\rho + (1 - \beta^2)^{3/2}}
\label{valorcriticokappa=0}
\end{equation}
and the numerator of $l_1$ in (\ref{coeficiente1}) is given by
{\small
\begin{eqnarray}\label{l1kappa=0}
G_2(\beta, \alpha, \rho) = -2 \alpha^4 \beta^{22} + 2 \alpha^4
\beta^{20} (8 + 7 \rho ((1 - \beta^2)^{1/2} + 3 \rho)) + 3 (-2+
\nonumber \\ \rho (-9(1 - \beta^2)^{1/2} + \rho(-15-10(1 -
\beta^2)^{1/2} \rho + 3 (1 - \beta^2)^{1/2} \rho^3 +\rho^4))) -
\nonumber \\ 2 \beta^{18} (-5+ \alpha^2 (1 + \rho ( (1 -
\beta^2)^{1/2} + 5 \rho)) + \alpha^4 (28 + \rho (50(1 -
\beta^2)^{1/2} + \nonumber \\ 7 \rho (22 + 5 \rho ((1 -
\beta^2)^{1/2} + \rho))))) + \beta^{16} (- 86+ \alpha^2 (16 + \rho
(27(1 - \beta^2)^{1/2} + \nonumber \\ 2 \rho (52 + 9 (1 -
\beta^2)^{1/2} \rho + 7 \rho^2))) + 2 \alpha^4 (56 + \rho (153 (1
- \beta^2)^{1/2} + \nonumber \\ 7 \rho (69 +\rho (33(1 -
\beta^2)^{1/2} + \rho (35 + 3(1 - \beta^2)^{1/2} \rho +
\rho^2)))))) - \nonumber \\ \beta^{14} (-41(8 + (1 -
\beta^2)^{1/2} \rho) + 2 \alpha^2 (28 + \rho (64 (1 -
\beta^2)^{1/2} + \rho (215 + \nonumber \\ 2 \rho (40(1 -
\beta^2)^{1/2} + \rho (34 + (1 - \beta^2)^{1/2} \rho))))+ \alpha^2
(70 + \nonumber \\ \rho (260 (1 -\beta^2)^{1/2}+ \rho (\rho (630
(1 - \beta^2)^{1/2} + 840 + \rho (700 + \nonumber \\ \rho (140 (1
- \beta^2)^{1/2} + \rho (56 + (1 - \beta^2)^{1/2} \rho))))))))) +
\beta^{12} (-13(56 + \nonumber \\ 21(1 - \beta^2)^{1/2} \rho + 4
\rho^2) + \alpha^2 (112 + \rho (319(1 - \beta^2)^{1/2} + \rho (976
+ \nonumber \\ \rho (552(1 - \beta^2)^{1/2} + \rho (494 + 51 (1 -
\beta^2)^{1/2} \rho + 8 \rho^2))))) + 2 \alpha^4 (56 + \nonumber
\\ \rho (265 (1 - \beta^2)^{1/2} + \rho(875 + \rho (910 (1 -
\beta^2)^{1/2} + \rho (1050 + \nonumber \\ \rho (350 (1 -
\beta^2)^{1/2}+ \rho (154 + 11(1 - \beta^2)^{1/2} \rho + \rho^2))))))))- \nonumber \\
\beta^{10} (-259(4 + 3(1 - \beta^2)^{1/2} \rho) + \rho^2 (-305 +
9(1 - \beta^2)^{1/2} \rho) + \alpha^2 (140 + \nonumber \\ \rho
(480(1 - \beta^2)^{1/2} + \rho (1376 + \rho (1011 (1 -
\beta^2)^{1/2} + 2 \rho (458 + \nonumber \\ 93(1 - \beta^2)^{1/2}
\rho + 23 \rho^2))))) + 2 \alpha^4 (28 +\rho (162 (1 -
\beta^2)^{1/2} + \rho (546 + \nonumber \\ \rho (735 (1 -
\beta^2)^{1/2} + \rho (875 + \rho (420 (1 - \beta^2)^{1/2} + \rho
(196 + \rho + \\ 3 \rho (9(1 - \beta^2)^{1/2}))))))))) +
\beta^8 (-245 (4 + 5 (1 - \beta^2)^{1/2} \rho) + \rho^2 (-745 + \nonumber \\
29 (1 - \beta^2)^{1/2} \rho + 27 \rho^2) + \alpha^2 (112 + \rho
(457 (1 - \beta^2)^{1/2} + \rho (1264 + \nonumber \\ \rho (1121 (1
- \beta^2)^{1/2} + 2 \rho (502 + \rho (143 (1 - \beta^2)^{1/2}+ 2
\rho (20 + (1 - \nonumber \\ \beta^2)^{1/2} \rho))))))) + 2
\alpha^4 (8 + \rho (55(1 - \beta^2)^{1/2} + \rho (189 + \rho (315
(1 - \beta^2)^{1/2} + \nonumber \\ \rho (385 + \rho (245 (1 -
\beta^2)^{1/2} + \rho (119 + 25(1 - \beta^2)^{1/2} \rho + 3
\rho^2)))))))) - \nonumber \\ \beta^6 (-77(8 + 15(1 -
\beta^2)^{1/2} \rho) + \rho^2 (-970+ 3(1 - \beta^2)^{1/2} \rho+ \nonumber \\
112 \rho^2) + \alpha^2 (56 + \rho (272 (1 - \beta^2)^{1/2} + \rho
(750 +\rho (795(1 - \beta^2)^{1/2} + \nonumber \\ \rho (710 + \rho
(240 (1 - \beta^2)^{1/2} + \rho (68 + 5(1 - \beta^2)^{1/2}
\rho))))))) + 2 \alpha^4 (1 + \nonumber \\ \rho (8(1 -
\beta^2)^{1/2} +\rho (28 + \rho (56(1 - \beta^2)^{1/2} + \rho (70
+\rho (56(1 - \beta^2)^{1/2} + \nonumber
\\ \rho (28 + 8(1 - \beta^2)^{1/2} \rho + \rho^2)))))))) + \beta^4 (-248 +
\rho (-651 (1 - \beta^2)^{1/2} + \nonumber \\ \rho (-710 + \rho
(-75(1 - \beta^2)^{1/2}+ \rho (143 +27 (1 - \beta^2)^{1/2}
\rho))))+ \alpha^2 (16 + \nonumber \\ \rho (93(1 - \beta^2)^{1/2}
+ \rho (264 +\rho (345(1 - \beta^2)^{1/2} + 135 (1 -
\beta^2)^{1/2} \rho^2 + \nonumber \\ 3(1 - \beta^2)^{1/2} \rho^4 +
40 \rho (8 + \rho^2)))))) - \beta^2 (-58+ \rho (-203(1 -
\beta^2)^{1/2} + \nonumber \\ \rho (-277 + \rho (-88(1 -
\beta^2)^{1/2} + \rho (58 + 47(1 - \beta^2)^{1/2} \rho + 9
\rho^2)))) + \nonumber \\ 2 \alpha^2 (1 + \rho (7 (1 -
\beta^2)^{1/2} + \rho (21 + \rho (35(1 - \beta^2)^{1/2} + \rho (35
+ \nonumber \\ \rho (21 (1 - \beta^2)^{1/2} + \rho (7 +(1 -
\beta^2)^{1/2} \rho)))))))) \nonumber.
\end{eqnarray}
}
If $G_2$ is different from zero then the family of HGS
differential equations (\ref{sistemafinal}) has a transversal Hopf
point at $P_0$ for $\varepsilon = \varepsilon_c$.

\label{teoremakappa=0}
\end{cor}

\noindent {\bf Proof.} The proof is obtained by substituting
$\kappa = 0$ into Eqs. (\ref{P0}), (\ref{valorcritico}) and
(\ref{coeficiente1}). The long expression above, being a challenge
to hand calculation, has been performed with Computer Algebra. In
the site \cite{mello} has been posted the main steps of the long
calculations  involved in this substitution. This has been done in
the form of a {\it notebook}  for MATHEMATICA 5 \cite{math}. A
sufficient condition for being a Hopf point is that the first
Lyapunov coefficient $l_1 \neq 0$, since the transversality
condition is satisfied by Proposition \ref{lematransv}. But from
(\ref{l1kappa=0}) it is equivalent to $G_2 \neq 0$.
\begin{flushright}
$\blacksquare$
\end{flushright}

\begin{figure}[!h]
\centerline{
\includegraphics[width=6cm]{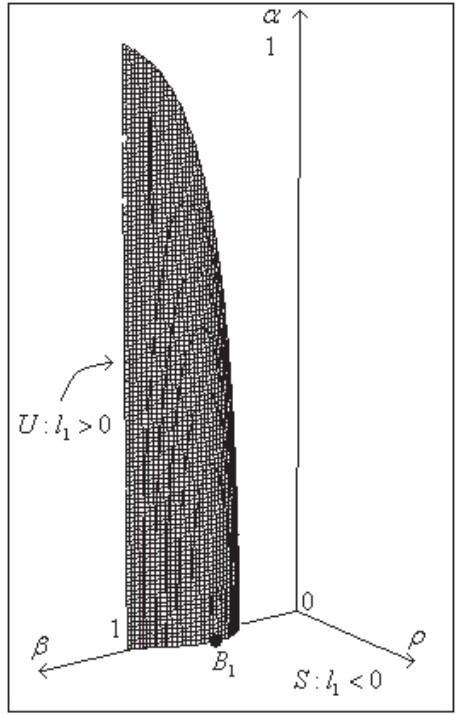}}

\caption{{\small Signs of the first Lyapunov coefficient for
$\kappa = 0$}.}

\label{sinalL1kappa0}
\end{figure}

Equation (\ref{l1kappa=0}) gives an expression to determine the
sign of the first Lyapunov coefficient (\ref{coeficiente1}) for
the case $\kappa = 0$. The graph $G_2(\beta, \alpha,\rho) = 0$ is
illustrated in Fig. \ref{sinalL1kappa0}, where the signs of the
first Lyapunov coefficient are also represented. The surface $l_1
= 0$ divides the hypersurface of critical parameters $\varepsilon
= \varepsilon_c$ into two connected components denoted by $S$ and
$U$ where $l_1 < 0$ and $l_1 > 0$ respectively. At point $B_1$ the
$\rho$ coordinate is $0.0478$. See Fig. \ref{sinalL1kappa0}.

The following theorem summarizes the results in this subsection.

\begin{teo}
Consider the case $\kappa = 0$. If $(\beta, \alpha, \rho) \in S
\cup U$ then the family of differential equations
(\ref{sistemafinal}) has a transversal Hopf point at $P_0$ for
$\varepsilon = \varepsilon_c$. If $(\beta, \alpha, \rho) \in S$
then the Hopf point at $P_0$ for $\varepsilon = \varepsilon_c$ is
asymptotically stable and for each $\varepsilon < \varepsilon_c$,
but close to $\varepsilon_c$, there exists a stable periodic orbit
near the unstable equilibrium point $P_0$. If $(\beta, \alpha,
\rho) \in U$ then the Hopf point at $P_0$ for $\varepsilon =
\varepsilon_c$ is unstable and for each $\varepsilon
> \varepsilon_c$, but close to $\varepsilon_c$, there exists an
unstable periodic orbit near the asymptotically stable equilibrium
point $P_0$. See Fig \ref{sinalL1kappa0}.

\label{teoremahopf1}
\end{teo}

\section{Concluding comments}\label{conclusion}

In this paper the original stability analysis due to Maxwell and
Vyshnegradskii of the Watt Centrifugal Governor System ---WGS---
has been extended to the Hexagonal Governor System ---HGS--- where
a more general force, due to the spring, acting on the sliding
sleeve of the governor has been considered.

In Theorem \ref{teoestabilidade} we have extended the stability
results presented in Pontryagin \cite{pon} to include this more general system. See \cite{smb1} for
another possible extension.

Concerning the bifurcations of the HGS, this paper deals with the
codimension one Hopf bifurcations in the Hexagonal governor
differential equations. The general expression for the first
Lyapunov coefficient at the Hopf point has been obtained in
Theorem \ref{teoremageral}. More concrete consequences of this
calculation have been synthesized in Theorems \ref{teoremahopf}
and \ref{teoremahopf1}. These results give sufficient conditions
for the stability of the points on the Hopf hypersurface and of
the periodic orbit that bifurcates from the Hopf point for the
Hexagonal governor differential equations (\ref{sistemafinal}) in
two particular cases easier to visualize with the help of
numerical plotting. See Figs. \ref{sinalL1rho0},
\ref{sinalL1kappa0} and the site \cite{mello}.

\vspace{0.2cm} \noindent {\bf Acknowledgement}: The first and
second authors developed this work under the project CNPq Grant
473824/04-3. The first author is fellow of CNPq and takes part in
the project CNPq PADCT 620029/2004-8.

\end{document}